\documentclass[a4paper,12pt]{amsproc}

\newtheorem{thm}{Theorem}
\newtheorem{lem}{Lemma}

\theoremstyle{remark}
\newtheorem{rk}{Remark}

\begin{document}

\title{ON THE MULTIPLICITY OF THE HYPERELLIPTIC INTEGRALS}
\author{Claire MOURA}
\address{Laboratoire de Math\'ematiques E. Picard,
UFR MIG, Universit\'e Paul Sabatier,
118 route de Narbonne,
31062 TOULOUSE cedex 4, FRANCE}

\email{moura@picard.ups-tlse.fr}

\begin{abstract}
Let $I(t)= \oint_{\delta(t)} \omega$ be an Abelian integral,
where $H=y^2-x^{n+1}+P(x)$ is a hyperelliptic polynomial of Morse type,
$\delta(t)$ a horizontal family of cycles in the curves $\{H=t\}$, and  
  $\omega$ a polynomial 1-form in the variables $x$ and $y$. We provide an
  upper bound on the multiplicity of $I(t)$, away from the critical values of
  $H$. Namely: $ord\ I(t) \leq n-1+\frac{n(n-1)}{2}$ if $\deg \omega <\deg
  H=n+1$.
The reasoning goes as follows: we consider the analytic curve
parameterized by
the integrals along $\delta(t)$ of the $n$ ``Petrov''
forms of $H$ (polynomial 1-forms that freely generate the module of
relative cohomology of $H$), and interpret the multiplicity of $I(t)$ as the
order of contact of $\gamma(t)$ and a linear hyperplane of $\textbf C^
n$. Using the Picard-Fuchs system satisfied by $\gamma(t)$, we establish an
algebraic identity involving the wronskian determinant of the integrals of the
original form $\omega$ along a basis of the homology of the generic fiber of
$H$. The latter wronskian is analyzed through this identity, which yields
 the estimate on the multiplicity of $I(t)$. Still, in some
cases, related to the geometry at infinity of the curves $\{H=t\} \subseteq
\textbf C^2$, the wronskian occurs to be zero identically. In this alternative
we show how to adapt the argument to a system of smaller rank, and get a
nontrivial wronskian. 
For a form $\omega$ of arbitrary degree, we are led to estimating the order of
contact between $\gamma(t)$ and a suitable algebraic hypersurface in $\textbf
C^{n+1}$. We observe that $ord\ I(t)$ grows like an affine function with
respect to $\deg \omega$.

\end{abstract}

\maketitle

\vspace{1cm}

 \section{Introduction}

Consider a complex bivariate polynomial $H(x,y) \in \textbf C[x,y]$. It is
well known that the polynomial mapping $H:\ \textbf C^2 \rightarrow \textbf C$
defines a locally trivial differentiable fibration over the complement of a finite subset of
$\textbf C$ (cf \cite{D}). Under some restrictions on the principal part of $H$, this set limits
to  $\mbox{crit}(H)$, the set of critical values of $H$ (see \cite{Br}, \cite{Gav1}). One can then consider the homology bundle
$\cup_t(H_1(\{H=t\}, \textbf Z)\rightarrow \textbf C\ \backslash\ \mbox{crit}(H)$, equipped with the
Gauss-Manin connection. Take a class $\delta(t)$ in the homology group
$H_1(\{H=t\}, \textbf Z)$ of a generic fiber of $H$. As the base of the
homology bundle is 1-dimensional, the connection is flat, and the parallel
transport of $\delta$ depends only on the homotopy class of the path in
$\textbf C \ \backslash \  \mbox{crit}(H)$. The transport of the homology
class $\delta(t)$ along a loop that
encircles a critical value of $H$ results in a nontrivial outcome, due to the
action of the monodromy on $\delta$. Let $\omega$ be a
polynomial 1-form in the variables $x$ and $y$. Its restriction on any 
fiber of $H$ is a closed form, therefore the integral of $\omega$ on a cycle
lying in a regular level curve of $H$ depends only on the homology class of this
cycle. Consider the complete Abelian integral $I(t)=\oint_{\delta(t)} \omega$. This function admits
an analytic extension in the complement of $\mbox{crit(H)}$. 
We refer to \cite{AGV}, \cite{Yak1} for a detailed survey.

Let $t_0 \in \textbf C$ be a 
regular value of $H$. One can ask for an estimate on the multiplicity of
$I(t)$ at $t_0$. The result is expected to depend on two parameters, namely  
the degree of $H$ and the degree of
$\omega$. The number of parameters can be reduced by looking first at the case
$\deg \omega< \deg H$.
This case is the most interesting regarding the connection
with the infinitesimal Hilbert 16th problem, that takes place in the real
setting. Assume $H$ has real coefficients. Consider the Hamiltonian
distribution $dH$, and the one-parameter perturbation $dH+ \epsilon \omega$ by an
arbitrary real polynomial 1-form $\omega$. The first order term in the Taylor
expansion at $\epsilon=0$ of the corresponding displacement function $d(t, \epsilon)$   
 is an integral $I(t)$ of $\omega$ along
an oval in a level curve of $H$. Under the assumption $\deg \omega< \deg H$,
it is proved by  Yu. Ilyashenko (cf with \cite
{Il2}) that: $I(t) \equiv 0$ if and only if $\omega$ is \textit{exact}, so that the
perturbation is still a Hamiltonian distribution. In the case when $I(t)
\not\equiv 0$, the multiplicity of $I(t)$ at a point $t_0$ provides an upper bound for the
cyclicity of the displacement function $d(t, \epsilon)$ at $(t_0, 0)$. Thus,
the vanishing of the Abelian integral is relevant to the number of limit cycles born by small perturbation
of the Hamiltonian distribution $dH$.

For polynomials $H$ with generic principal part ($H$ regular at infinity), an
answer on the order of $I(t)$ is given by P. Mardesic in \cite{Mar}: a step towards the multiplicity of the
Abelian integrals consists in measuring the multiplicity of their Wronskian
determinant, which is a globally univalued function on $\textbf {CP}^1$, hence rational, with poles at the critical values of $H$ and possibly at infinity.  

We focus on the Abelian integrals performed on level curves of
\textit{hyperelliptic} polynomials. We establish a relation between the
Wronskian and a polynomial that we build up from the Picard-Fuchs system. 
As the Picard-Fuchs system reflects the topology of the level sets of $H$,
our approach depends whether $\deg H$ is even
or odd, still our estimate on the multiplicity of $I(t)$  is always quadratic with respect to $\deg H$. We finally show that, for a fixed hyperelliptic Hamiltonian, the growth of the multiplicity of $I(t)$ is
linear with respect to $\deg \omega$.

\section{Preliminary observations}

We begin by recalling a result about flatness of solutions of a linear
differential system. Consider a system $dx= \Omega x$ of order $n$, whose
coefficient matrix $\Omega$ is meromorphic on $\textbf{CP}^1$. Denote by
$t_1, \ldots, t_s$ the poles of $\Omega$. Fix a point $t_0$, distinct from the
poles, and consider a solution $\gamma(t) \subseteq \textbf C^n$, analytic in a
neighbourhood of $t_0$. Take a linear hyperplane $\{h=
\sum_{i=1}^{n}c_ix_i=0\} \subseteq \textbf C^n$. If this hyperplane does not
contain the solution $\gamma$, then (cf \cite{Mou}):

\begin{thm}
$$ord_{t=t_0}(h \circ \gamma)(t) \leq n-1 + \frac{n(n-1)}{2}
\left(\sum_{i=1}^{s}(-ord_{t_i}\Omega)-2 \right)$$ where $ord_{t_i}\Omega$ is
the minimum order of the pole $t_i$ over the entries of $\Omega$.
\end{thm}

\bigskip

We give here a simplified algorithm of the proof: write the system in the affine
chart $t$ in the form
$\dot x=\frac{A(t)}{P(t)} x$,for a polynomial matrix $A$ and scalar polynomial $P$. Replace
the derivation $\frac{\partial}{\partial t}$ by $D=P(t)
\frac{\partial}{\partial t}$. Then the curve $\gamma$ satisfies: $D \gamma=A \gamma$. Due to the linearity, we can write
$y(t)= (h \circ \gamma)(t)$ as the product of the row matrix $q_0=(c_1, \ldots, c_n)$
by the column matrix $\gamma$: $y(t)= q_0 \cdot \gamma(t)$. The successive derivatives of $y$ with
respect to $D$ can be written in a similar way: $D^ky(t)= q_k(t) \cdot \gamma(t)$,
where the row vectors $q_k$ have polynomial coefficients and are constructed
inductively by: $q_{k+1}=Dq_k+q_kA$. We observe that the sequence of $\textbf C(t)$-vector
spaces $V_k \subseteq \textbf C(t)^n$ spanned by the
vectors $q_0,q_1, \dots, q_k$, is strictly increasing (before stabilizing), hence we may extract from the matrix $\Sigma$ with rows $q_0, \ldots, q_{n-1}$, a nondegenerate minor $\Delta$ of rank $l \leq n$, such that any vector $q_k$ decomposes according to the Cramer rule:
$$q_k(t)= \sum_{i=0}^{l-1} \frac{p_{ik}(t)}{\Delta(t)} q_i(t),\quad k\geq l$$ with \textit{polynomial} coefficients $p_{ik}(t)$.
This shows that the function $y$ is a solution of an infinite sequence of
linear differential equations of the form:

\begin{equation}
\Delta \cdot D^ky=\sum_{i=0}^{l-1}p_{ik}(t)D^iy,\quad k\geq l
\end{equation}

Then, by deriving in an appropriate way each of these relations, one arrives
at the key-assertion: 
$$ord_{t_0}y \leq l-1+ ord_{t_0} \Delta$$
Thus, the flatness of a particular solution is correlated to the multiplicity of a polynomial constructed from the system.

One can derive an analytic version of this assertion by complementing the solution $\gamma$ by $n-1$ vector-solutions $\Gamma_2, \ldots, \Gamma_{n-1}$, so as to obtain a
fundamental matrix in a simply
connected domain around $t_0$. In particular: $\det(\gamma, \Gamma_2,\ldots, \Gamma_n)(t)$ does not vanish in this domain.

We restrict these solutions on the hyperplane $\{h=0\}$ and
set: $y_1=y=(h \circ \gamma)$, $y_i=(h \circ \Gamma_i), \ i=2,
\ldots,n$.  Let $l \leq n$ be the maximum number of
independent functions among $y_1, \ldots, y_n$.
  Their Wronskian determinant $W=W(y_1, \ldots, y_l)$ is analytic and does not
 vanish   identically around $t_0$. Expand $W$ with respect to any of its
columns, it follows that:
$$ord_{t_0}W \geq min_{k=0, \ldots, l-1}
\{ord_{t_0}y_i^{(k)}\}+ord_{t_0} D_k$$ where $D_k$ is the minor corresponding to the
element $y_i^{(k)}$. Hence: $ord_{t_0}W \geq ord _{t_0}y_i-(n-1)$, for any $i=1,
\ldots, n$. So: 
$$ord_{t_0}(h\circ \gamma) \leq n-1 + ord_{t_0}W$$ 

Naturally, the order of vanishing of $W$ does not depend on the particular
choice of fundamental system.  Besides, one arrives at the same conclusion by
forming the Wronskian determinant $W_D(y_1, \ldots, y_l)$ with respect to the
derivation $D=P(t) \frac{\partial}{\partial t}$, since  $W_D=P^{l(l-1)/2}\cdot
W$, and $D$ is not singular at $t_0$ (meaning that
 $P(t_0) \neq 0$).

Note that, like in the algebraic situation, one can interpret  $W_D$ as the principal coefficient of a linear differential equation
of order $l\leq n$ satisfied by $y_1, \ldots, y_l$:

\begin{equation}
 W_D (y_1, \ldots, y_l)D^ly+a_{l-1}(t) D^{l-1}y+ \ldots + a_0(t)y=0
\end{equation}
with coefficients $a_i(t)$ analytic in a neighbourhood of ${t_0}$. The method
 leading to such an equation is standard.   
For any linear combination $y$ of $y_1, \ldots, y_l$, the following Wronskian
determinant of size $l+1$ is zero identically: 

$$\begin{array}{|cccc|}
y_1 & \ldots & y_l & y\\
Dy_1 & \ldots & Dy_l & Dy\\
\vdots & \vdots & \vdots & \vdots \\
D^{l-1}y_1 & \ldots & D^{l-1} y_l& D^{l-1} y\\
D^{l}y_1 & \ldots & D^{l}y_l & D^{l}y
\end{array}$$

\bigskip

Expanding this determinant with respect to its last column gives the
equation. This is the analytic analogue of the $lth$ order equation in the
sequence 
(1). Both of them admit the same solutions.

\bigskip

Suppose that $\Sigma$ is a nondegenerate matrix, that is, its determinant
$\Delta(t)$ is not the null polynomial. Let $\mathcal P$ be a fundamental matrix
of solutions of the system $Dx=Ax$. 
From the construction, we obtain immediately the following matrix relation:

\bigskip

\begin{lem}
$\Sigma \cdot \mathcal P =  \mathcal W_D(y_1, \ldots, y_n)$, where $\mathcal
W_D(y_1, \ldots, y_n)$ is the Wronski matrix of $y_1, \ldots, y_n$, computed with the derivation $D$. 
\end{lem}

\bigskip

\begin{rk}
The matrix $\Sigma$ defines a meromorphic gauge equivalence between the
original system and the companion system of the equation (2).
\end{rk}

\bigskip

This yields the relation between determinants:
\begin{equation}
\det \Sigma \cdot \det \mathcal P= W_D(y_1, \ldots, y_n)= P(t)^
{\frac{n(n-1)}{2}} \cdot W(y_1, \ldots, y_n)
\end{equation}

$W(y_1, \ldots, y_n)$ being
the usual Wronskian $W_{\frac{\partial}{\partial t}}(y_1, \ldots, y_n)$.

\bigskip

Now, at the non-singular point $t_0$ of the system, both $P$ and
$\det\mathcal P$ are nonzero, so that the order at  $t_0$ of the analytic function $W(y_1, \ldots, y_n)$  is
exactly the order at  $t_0$ of the polynomial determinant $ \det \Sigma$. 

\section{Multiplicity of the integrals}

\subsection{Petrov forms and Picard-Fuchs system}

Consider a bivariate hyperelliptic polynomial $H \in \textbf C [x,y]$, $H=y^2-x^{n+1}+ \overline H(x)$, with $\deg \overline H=n-1$.
Hyperelliptic polynomials are examples of semi quasi-homogeneous
polynomials. Recall that a polynomial $H$ is said to be semi quasi-homogeneous if the
following holds: the variables $x$ and $y$ being endowed with weights $w_x$
and $w_y$ (so that a monomial $x^\alpha
y^\beta$ has weighted degree $\alpha w_x+\beta w_y$), $H$ decomposes as a sum
$H^*+ \overline H$, and the highest weighted-degree part $H^*$ possesses an
isolated singularity at the origin. Moreover, a semi quasi-homogeneous
polynomial has only isolated singularities.
In the sequel, the notation $``\deg''$ will
stand for the weighted degree. For $H$ hyperelliptic, $\deg H=n+1$, with
$w_x=1$, $w_y=(n+1)/2$.  The weighted degree
extends to polynomial 1-forms: for $\omega=P(x,y) dx+Q(x,y)dy$, $\deg
\omega$ is the maximum $max(\deg P+w_x, \deg Q+w_y)$.
The symbol  $\Lambda^k$ will designate the $\textbf C [x,y]$-module of
polynomial k-forms on $\textbf C^2$.

Consider the quotient
 $\mathcal P_H= \frac{\Lambda^1}{\Lambda^0dH+d\Lambda^0}$. It is a module over
 the ring of polynomials in one indeterminate. Note that the integral of a
 1-form in $\Lambda^1$ depends only on its class in $\mathcal P_H$. In
 addition, working in the Petrov module of $H$ enables to exhibit a finite
 number of privileged 1-forms, that we will call the Petrov forms: indeed, $\mathcal P_H$ is freely generated by the monomial 1-forms $\omega_1=ydx$, $\omega_2=xydx, \ldots$, $\omega_n=x^{n-1}ydx$. 
 Moreover, the class of any 1-form in $\mathcal P_H$
decomposes as a sum: $p_1(t) \omega_1+ \ldots+ p_n(t) \omega_n$, with the following estimates on the degrees of the polynomials $p_i$:
$$\deg p_i \leq \frac{\deg \omega - \deg
  \omega_i}{\deg H}$$
These assertions belong to a general theorem due to L. Gavrilov (\cite{Gav2}),
where the Petrov module of any semi quasi-homogeneous polynomial $H$ is
described. The number of Petrov forms is the global Milnor number of $H$.

\bigskip

Consider a hyperelliptic integral  $\oint_{\delta(t)} \omega$, in a
neighbourhood of a regular value $t_0$ of the Hamiltonian $H$.
It becomes natural to consider the germ of analytic curve $\gamma(t)=\left(\oint_{\delta} \omega_1, \oint_{\delta} \omega_2,
\ldots, \oint_{\delta} \omega_n \right)$ parameterized by the integrals of the Petrov
forms. We shall start with forms of small degree, that is $\deg \omega \leq n$, and  study the behaviour of the multiplicity of the integral with
respect to  $n$. This restriction on the degree implies that $\omega$ is a linear
combination of the Petrov forms, with \textit{constant} coefficients:
$\omega=\sum_{i=1}^n c_i \omega_i$, $c_i \in \bf C$. Therefore,  the question
amounts to estimating the order of contact at the point $t=t_0$ of the curve
$\gamma(t)$ and of  the linear hyperplane
$\{\sum_{i=0}^{n}c_i x_i=0\}$ whose coefficients
are prescribed by the decomposition of $\omega$. 

\bigskip

In order to apply the argument presented in Section 2, we have to interpret
$\gamma$ as a solution of a linear differential system. We recall the
procedure described by S. Yakovenko in \cite[Lecture 2]{Yak1}. 
For any $i=1, \ldots, n$, divide the 2-forms $Hd\omega_i$ by $dH$. Then apply
the Gelfand-Leray formula and decompose the Gelfand-Leray residue in the
Petrov module of $H$. For a hyperelliptic Hamiltonian, it is clear that the $\textbf C$-vector space of relative 2-forms $\frac{\Lambda^2}{dH \wedge \Lambda^1}$
is spanned by the differentials of the Petrov forms $d\omega_1, \ldots, d\omega_n$. Whence the decomposition:

\begin{equation}
 H \cdot d\omega_i= dH \wedge \eta_i+\sum_{j=1}^{n}a_{ij}d\omega_j, \
a_{ij}\in \bf C,\  \eta_i \in \Lambda^1.
\end{equation}

Now, by the Gelfand-Leray formula,
\begin{equation}
 t\frac{d}{dt}\oint_\delta
\omega_i-\sum_{j=1}^{n} a_{ij}\frac{d}{dt} \oint_\delta \omega_j= \oint_\delta
\eta_i
\end{equation}
and this relation does not depend on the cycle of integration.

In order to obtain the system, one decomposes the residues $\eta_i$ in $\mathcal
P_H$. So, an estimate on the degree of $\eta_i$ is required, which can be
quite cumbersome when starting from a general Hamiltonian. Yet, in the
hyperelliptic case, (4) is completely explicit (cf with \cite{NY}), and one
sees immediately that: for any $i$, $\eta_i= \sum_{j=1}^{n}b_{ij} \omega_j$,
$b_{ij} \in \textbf C$.
Then (5) appears as the expanded form of the linear system:
$$(tE-A) \dot x = B x$$
where $A=\left ( a_{ij} \right )$ and $B=\left ( b_{ij} \right )$ are constant matrices, and $E$ is
the $(n\times n)$ Identity matrix.

We can write it as well as a system with rational
coefficients: $\dot x=\frac{C(t)}{P(t)} x$, where the polynomial matrix $C(t)$, obtained as $C(t)=\mbox{Ad}(tE-A) \cdot B$, has degree $n-1$,  and the scalar polynomial $P(t)= \det (tE-A)$ has degree $n$. The
polynomial $P$ can be explicited: evaluation of the relation (4) at a critical
point $(x_*,y_*)$ of $H$ shows that the corresponding critical value $t_*=H
(x_*,y_*)$ is an eigenvalue of $A$.
If the critical values of $H$ are assumed pairwise distinct, then:
$P(t)= (t-t_1) \ldots (t-t_n)$. Thus,
the singular points of the system are the critical values of $H$ and the point
at infinity. All of them are Fuchsian.

\bigskip

We can apply Theorem 1 and get the following estimate on the multiplicity at a
zero of the integral: $ord_{t_0} \oint_\delta \omega \leq n-1+
\frac{n(n-1)^2}{2}$. We are going to show how to improve this bound, applying
(3).

\subsection{Main result}

We now formulate the theorem. We impose an additional requirement on the
hyperelliptic Hamiltonian: $H$ has to be of Morse type, that is, with nondegenerate
critical points as well as distinct critical values.

\begin{thm}
Let $H$ be a hyperelliptic polynomial $H=y^2-x^{n+1}+ \overline H(x)$,
where $\overline H$ is a
polynomial of degree $n-1$, of Morse type. Let $\{\omega_1, \ldots,
\omega_n\}$ be the set of monomial Petrov forms associated to $H$, and let
$\omega= \sum_{i=1}^{n}c_i \omega_i$, $c_i \in \textbf C$ be an arbitrary
linear 
combination. Consider the Abelian integral of $\omega$ along a horizontal
section $\delta(t)$ of the homology bundle.  
 Let $t_0$ be a regular value of $H$. If
$\oint_{\delta(t)} \omega \not \equiv 0$, then:

$$ord_{t_0}\left( \oint_{\delta(t)} \omega \right)\leq n-1+\frac{n(n-1)}{2}$$
\end{thm}

\bigskip

The rest of this section is devoted to the proof of Theorem 2. 

Recall that a fundamental matrix
 of the Picard-Fuchs system is obtained by integrating some $n$ suitable polynomial 1-forms $\Omega_1, \ldots, \Omega_n$ over any  basis of the homology
groups $H_1(\{H=t\}, \textbf Z)$:
$$\mathcal P=
\left (
\begin{array}{ccc}
\oint_{\delta_1}  \Omega_1 & \cdots & \oint_{\delta_n} \Omega_1\\
\vdots & \vdots & \vdots \\
\oint_{\delta_1}\Omega_n & \cdots & \oint_{\delta_n} \Omega_n
\end{array}
\right )$$

As explained in \cite[Lecture 2]{Yak1}, the determinant of a fundamental
 system of solutions has to be a polynomial, divisible by $(t-t_1) \ldots
 (t-t_n)$. Its actual degree depends on the choice of the integrands.
It is shown by L. Gavrilov in \cite{Gav2} and D. Novikov in \cite{Nov} that
 one can plug the Petrov forms into the period matrix $\mathcal P$ and get
 $\det \mathcal P= c \cdot (t-t_1) \ldots (t-t_n)=P(t)$, with a
\textit{nonzero}
constant $c$.

Next, we form the vectors $q_0, \ldots, q_{n-1}$, $q_0=(c_1, \ldots, c_n)$,
$q_{k+1}=Dq_k+q_k C(t)$, $D=P(t) \frac{\partial}{\partial t}$, and collect
them in a matrix $\Sigma$. Lemma 1 says:

$$\Sigma \cdot \mathcal P= \mathcal W_D \left(\oint_{\delta_1} \omega, \ldots,
\oint_{\delta_n} \omega \right)$$

which gives:

\begin{equation}
\Delta \cdot (t-t_1) \ldots (t-t_n)= P^{n(n-1)/2} \cdot W \left(\oint_{\delta_1} \omega,
\ldots,\oint_{\delta_n} \omega\right)
\end{equation}

\bigskip

The disadvantage of this formula is that it does not resist possible degeneracy of the
matrix $\Sigma$: the determinant $\Delta= \det \Sigma$ is a polynomial in the variable $t$, whose coefficients are homogeneous polynomials with respect to the components of $q_0$.
$$ \Delta_{q_0}(t)=P_0(q_0)+P_1(q_0) \cdot t+ \ldots +P_D(q_0) \cdot t^D$$
where $D$ is the maximum possible degree for $\Delta$, achieved for generic $q_0$.
One cannot a priori guarantee that the algebraic subset 
$$S=\{q_0 \in \textbf C^n:\ P_i(q_0)=0\ , i=0, \ldots,D\}$$ is reduced to zero.

We are now going to analyze on what conditions the equality (6) makes
sense. It involves 
the geometry at infinity of the hyperelliptic affine curves
$\{H=t\} \subseteq \textbf C^2$. Suppose that for the form $\omega=\sum_i ^{}
c_i \omega_i$, the Wronskian $W(\oint_{\delta_1} \omega,
\ldots,\oint_{\delta_n} \omega)$ vanishes identically as a function of
$t$. This means that one can find
a cycle $\sigma$, complex  combination of $\delta_1,
\ldots, \delta_n$, such that the integral $\oint_\sigma \omega$ is identically
zero. 

\bigskip

\begin{lem}
If an Abelian integral $\oint_{\sigma(t)} \omega$ is zero identically, then
$\sigma(t)$ belongs to
the kernel of the intersection form on $H_1(\{H=t\}, \textbf C)$.
\end{lem}

\begin{proof}
Assume on the contrary, that $\sigma$ has a nonzero intersection number with a
cycle from $H_1(\{H=t\}, \textbf C)$, while  $\oint_\sigma \omega$ vanishes identically.
The semi quasi-homogeneity property implies that $H$ defines a trivial
fibration at infinity, and the homology of a
regular fiber can be generated by a basis of $n$ vanishing cycles $v_i$ ( $v_i$
contracts to a point when t approaches $t_i$). Necessarily, $\sigma$
intersects one of the vanishing cycles  - say $v_1$: $(\sigma, v_1) \neq 0$. 
Continue analytically the integral $\oint_\sigma \omega$ along a loop around $t_1$: the Picard-Lefschetz formula states that the monodromy changes $\sigma$ into $\sigma - (\sigma, v_1) v_1$.  On the level of the integral: $0=0- (\sigma, v_1) \oint_{v_1} \omega$, hence $\oint_{v_1} \omega$ is zero.

 Moreover, the assumptions of hyperellipticity and Morse type imply that one can produce a
basis  $\{v_1, \ldots, v_n\}$ of vanishing cycles in which any two consecutive
cycles intersect: $(v_i,v_{i+1})=\pm 1$ (cf \cite{Il1}).
One proceeds inductively with the rest of the critical values $t_2, \ldots, t_n$:
 every integral of $\omega$ along $v_2, \ldots, v_n$ vanishes identically as well. This means that the \textit{restrictions} of $\omega$ on any
generic fiber of $H$ are exact.
 We can now apply L. Gavrilov's result \cite[Theorem 1.2]{Gav2} and deduce the global statement: the form $\omega$ has to be
 \textit{exact}, hence 
zero in the Petrov module. But this is clearly impossible since $\omega$ is a
combination of $\textbf C[t]$-independent forms.
\end{proof} 

Consequently, if the integral of $\omega$ vanishes along $\sigma$, this means
that $\sigma$ becomes homologous to zero when the affine fiber is embedded in its normalization.
That is, $\sigma$ lies in the kernel of the morphism $i_*$: $H_1(\Gamma, \textbf Z) \rightarrow
H_1(\overline \Gamma, \textbf Z)$, denoting  the affine curve $\{H=t\}
\in \textbf C^2$ by $\Gamma$ and its normalized curve in $\textbf {CP}^2$ by
$\overline \Gamma$.

\begin{lem}
If $n$ is even, then $S$ is limited to $\{0\}$.
\end{lem}

\begin{proof}
The projection
$\Pi:$ $\overline \Gamma \rightarrow \textbf{CP}^1$, $(x,y) \mapsto x$ is a
double ramified covering of $\textbf{CP}^1$. 
From the affine equation $y^2=\Pi_{i=1}^{n +1}(x-x_i(t))$, one finds $n+1$
ramification points of $\Pi$ in the complex plane, hence the total number of
ramification points of $\Pi$ is $n+1$ or $n+2$. On the other hand, using
the Riemann-Hurwitz formula, the number of ramification points of $\Pi$ is
 $2g_{\overline \Gamma}+2$. This shows that the genus $g_{\overline
  \Gamma}$ is equal to $[n/2]$. 
Therefore, if $n$ is even, the homology group of $\overline \Gamma$ has rank
$2g_{\overline
  \Gamma}=n$, and $i_*$ is an isomorphism. In this case, there is a single point at
infinity on $\overline \Gamma$ above $\infty \in \textbf {CP}^1$ and $\sigma$
is zero in the homology of the \textit{affine} level curve $\{H=t\}$. This
means that no relation can occur between the integrals $\oint_{\delta_1}
\omega, \ldots, \oint_{\delta_n} \omega$, unless $\omega$ is zero. 
\end{proof}

\bigskip

For even $n$, we can carry out the analysis further. The Wronskian $W=
W \left(\oint_{\delta_1} \omega,
\ldots,\oint_{\delta_n} \omega\right)$ is a rational function, so the sum of
its orders at all points of  $\textbf{CP}^1$ equals $0$. As a consequence, the order at one of its zeros can be deduced from the order at its poles and at the point at infinity:
$$ord_{t_0}W \leq -ord_\infty W -\sum_{i=1}^n ord_{t_i} W$$

From (6), we get: 

\begin{equation}
ord_\infty W =  ord_\infty \Delta -n + \frac{n^2(n-1)}{2}
\end{equation}

One gets easily an upper bound on $\deg \Delta$: from the inductive
construction the degree of each component of a vector $q_k$ is no larger than
$k(n-1)$, this yields: $\deg \Delta \leq \frac{n(n-1)^2}{2}$. In the right
hand side of (7), the cubic terms
cancel out each other, so that we get the estimate: 
 $ord_\infty W \geq \frac{n^2-3n}{2}$. This shows in particular that $\infty$
 is a zero of $W$.

As for the order of the $W$ on the finite singularities, we reproduce an
argument due to P. Mardesic (\cite{Mar}): the critical points of $H$ are Morse, which allows to
fix  the Jordan structure of the monodromy matrix at $t_i$, by choosing an
adapted basis of cycles. This imposes the structure of the integrals in a
neighbourhood of $t_i$: they are all analytic at $t_i$, except one of them that undergoes
ramification. The pole of the Wronskian at $t_i$ may only result from the derivation of this
integral. The estimate  follows automatically: $ord_{t_i} W \geq 2-n$. The contribution of the poles is $\sum_{i=1}^n ord_{t_i} W \geq 2n-n^2$. Therefore, we have obtained the following upper bound:
$$ord_{t_0}W \leq \sum_{t_0 \in \bf C \bf P^1, t_0
  \not= \infty,\ t_0 \not=t_i} ord_{t_0}W \leq \frac{n(n-1)}{2}$$
which proves Theorem 2 for even $n$.

\bigskip

We now return to the case of odd $n$.
The homology group $H_1(\overline \Gamma, \textbf Z)$ has rank
 $2g_{\overline \Gamma}=n-1$, and $\ker\ i_*$ is generated by one cycle that we will call $\delta_\infty$.
The forms $\omega$ that annihilate the Wronskian $W$ are those with zero integral along the cycle $\delta_\infty$. In order to describe this subspace of forms, we prove a dependence relation among Abelian integrals:

\begin{lem}
The complex vector space generated by the residues of the Petrov forms at infinity has dimension 2, that is:
$$\dim_\textbf C \left( \oint_{\delta_\infty} \omega_1, \ldots, \oint_{\delta_\infty} \omega_n \right)=2$$

\end{lem}

\begin{proof}
In order to estimate the residues $\rho(\omega_i)$ of the forms $\omega_i=x^{i-1}ydx$, $i=1, \ldots,
n-1$, at the point at infinity on the curve $y^2- x^{n+1}+ \overline H(x)-t=0$,
($[0:1:0] \in \textbf{CP}^2$), we pass to the chart $u=1/x$ and, expressing
$y$ as a function of $x$: $y=\pm (x^{n+1}+ \overline H(x)-t)^{1/2}$, we get
meromorphic 1-forms at $u=0$: 
$$\omega_i=  (1/u)^{i-1} \cdot  (1/u)^2 \cdot
(1/u)^{(n+1)/2} \cdot (1+R(u)+tu^{n+1})^{1/2} du$$ where $R$ is a
polynomial, $\deg R=n+1$, $R(0)=0$.
We have to compute the coefficient of $1/u$. The Taylor expansion of the
square root gives: for $i< \frac{n+1}{2}$, $\rho(\omega_i)$ is a constant with
respect to $t$, and for  $i \geq \frac{n+1}{2}$, $\rho(\omega_i)$ is a
polynomial of degree 1. This proves that the space of residues has dimension 2
over $\textbf C$, and is generated by any pair $\{\rho(\omega_i),
\rho(\omega_{i+(n+1)/2)})\}$, $i< \frac{n+1}{2}$.    
\end{proof}

\bigskip

\begin{rk}

In the case of a Hamiltonian of even degree (that is, for odd $n$), we detect
solutions of the Picard-Fuchs system that are included in a
hyperplane. Indeed, 
whenever a form $\omega=\sum_{i=1}^{n} c_i \omega_i$ has a zero residue at infinity,
then its coefficients define the equation of a hyperplane $\{h=\sum_{i=1}^{n}
c_i x_i=0\}$ that contains the integral curve 
 $\Gamma_1(t)=(\oint_{\delta_\infty}\omega_1, \ldots,
\oint_{\delta_\infty}\omega_n)$. This implies that the global monodromy of the
Picard-Fuchs 
system is reducible: extend $\Gamma_1$ to a fundamental system by adjoining solutions $\Gamma_2, \ldots, \Gamma_n$. Then, the $\textbf
C$-space spanned by the solutions $\Gamma_i$ such that $h \circ \Gamma_i(t) \equiv 0$
is invariant by the monodromy (see \cite[Lemma 1.3.4]{Bo}).

\end{rk}

\bigskip

It follows from Lemma 4 that the set $S$ coincides with the set of relations
between the residues of the Petrov forms at infinity. It is thus a codimension 2 linear subspace of $\textbf C^n$.
Therefore, we arrive at the conclusion of Theorem 2, for odd $n$ and $\omega
\not \in S$. In the remaining cases, that is for $\omega \in S$, the identity
(6) is useless, since both sides are $0$. 

\bigskip

We  aim at reconstructing the
identity (6), initiating the reasoning from a linear differential system of
size smaller than $n$. First, we know the exact number of independent integrals
among
 $\oint_{\delta_1}\omega, \ldots,\oint_{\delta_n}\omega$. 

\begin{lem}
If $\omega$ belongs to $S$, then:
$\dim_{\textbf C} \left(\oint_{\delta_1}\omega, \ldots,\oint_{\delta_n}\omega\right)=n-1$.
\end{lem}

\begin{proof}
The relations between these integrals constitute the space
$$\left \{(d_1, \ldots, d_n) \in
\textbf C^n:\  d_1 \oint_{\delta_1}\omega +\ldots+ d_n  \oint_{\delta_n}\omega \equiv
0 \right\}$$ From Lemma 2, any relation $(d_1, \ldots, d_n)$ must verify: $d_1 \delta_1+ \ldots + d_n \delta_n$ is a multiple of
$\delta_\infty$. 
This defines a 1-dimensional vector space.
\end{proof}

\bigskip

From now on, we work with an adapted basis of cycles in $H_1(\{H=t\},
\textbf Z)$, that includes $\delta_\infty$, and that we denote by
$(\delta_1, \ldots, \delta_{n-1}, \delta_\infty)$. Then it is clear that
$W(\oint_{\delta_1}\omega, \ldots,\oint_{\delta_{n-1}} \omega) \not \equiv 0$. We also
make several changes in the Petrov frame: recall that the matrix $A=(a_{ij})$ in (4)
describes the vectors $H d\omega_i$ via the correspondence that associates the $i$-th
canonical vector  $e_i \in \textbf C^n$ to the monomial $x^i$. We already
noticed that the matrix $A$ was diagonalizable.
 We assume $A$ diagonal (which
corresponds to combining linearly the forms $\omega_i$).
Moreover, we know that among the forms $\omega_i$ associated to an eigenbasis
of $A$, two of them will have independent residues at infinity. Up to
permutation, these are $\omega_{n-1}$ and $\omega_n$.
 Now,
after adding a scalar multiple of $\omega_n$ to $\omega_{n-1}$, we may assume
that 
$\oint_{\delta_\infty} \omega_{n-1}$ is a constant, while the residue of the
form $\omega_n$ is a polynomial of degree $1$ in the variable $t$. With respect to such a
basis $(\omega_1, \ldots, \omega_n)$, a nonzero off-diagonal entry, $a_{n,n-1}$,  may appear in $A$. Thus,
with our choice of basis of $\textbf C^n$, $A$ has the form: 

$$A=
\left (
\begin{array}{cccc}
a_{1,1} & & & \\
 & \ddots & & \\
 & & a_{n-1,n-1}& \\
 & & a_{n,n-1} &  a_{n,n} \\
\end{array}
\right )$$

Now, a  form
$\omega=\sum_{i=1}^{n}c_i \omega_i$ belongs to $S$ if and only if 
$$c_1 \oint_{\delta_\infty} \omega_1+ \ldots +c_n \oint_{\delta_\infty}
\omega_n \equiv 0$$ 
So that after a linear change of coordinates in $\textbf C^n$, we may write
$\omega$ as: $\omega=c_1 \widetilde \omega_1+ \ldots + c_{n-2} \widetilde
\omega_{n-2}$, with $\oint_{\delta_\infty}\widetilde \omega_i=0$, $i=1, \ldots, n-2$.
We set: $\widetilde \omega_{n-1}= \omega_{n-1}$ and $\widetilde \omega_{n}=
\omega_{n}$. Thus, the integral
$\oint_{\delta} \omega$ reads: $q_0 \cdot \overline{\gamma} (t)$, $q_0
\in \textbf C^{n-2}$, $\overline{\gamma} (t)=(\oint_\delta\widetilde \omega_{1}, \ldots, 
\oint_\delta \widetilde \omega_{n-2})$.

When passing to the Petrov frame $(\widetilde \omega_{i})$, $i=1, \ldots, n$,
the matrix $A$ is changed into $A'=P^{-1}AP$, and part of the
structure of the matrix $P=(p_{i,j})$ is known: $p_{n-1,n-1}=p_{n,n}=1$;
 $p_{i,n-1}=0$ for $i \neq n-1$, and
 $p_{i,n}=0$ for $i \neq n$, which implies that:
$a'_{i,n-1}=0$ for $i \neq n-1$ and $a'_{i,n}=0$ for $i \neq n$.

We now write the corresponding decomposition of the $n-1$ first 2-forms $Hd
\widetilde \omega_i$ and perform integration \textit{along the cycle}
$\delta_\infty$: 

$$ t \frac{d}{dt} \oint_{\delta\infty} \widetilde \omega_i -\sum_{j=1}^{n-1}
a'_{i,j}  \frac{d}{dt} \oint_{\delta\infty} \widetilde \omega_j = \oint
_{\delta\infty} \widetilde \eta_i, \quad i=1, \ldots, n-1$$
Since the residues of $\widetilde \omega_1, \ldots, \widetilde \omega_{n-1}$
are constants, these equalities entail: $\oint
_{\delta\infty} \widetilde \eta_i=0$, $i=1, \ldots, n-1$. Hence the Gelfand-Leray forms
$\widetilde \eta_i$ admit a decomposition with respect to $\widetilde \omega_1, \ldots, 
\widetilde \omega_{n-2}$ only: $\widetilde \eta_i= \sum_{j=1}^{n-2} b'_{i,j} \widetilde
\omega_j$, $i=1, \ldots, n-1$, $b'_{i,j} \in \textbf C$.
 Besides, from the
$n$-th equality
$$  t \frac{d}{dt} \oint_{\delta\infty} \widetilde \omega_n -
a'_{n,n}  \frac{d}{dt} \oint_{\delta\infty} \widetilde \omega_n = \oint
_{\delta\infty} \widetilde \eta_n$$
it follows that $\widetilde \eta_n$ has a nonzero component along $\widetilde \omega_n$.  
      
This provides information on the matrix $B'$ related to the new frame $(\widetilde \omega_i)$: $b'_{i,n-1}= b'_{i,n}=0$ for $i=1, \ldots, n-1$.

\bigskip

The curve $\gamma(t)=(\oint_\delta \widetilde \omega_1, \ldots, \oint_\delta
\widetilde \omega_n)$ is a solution of the linear system $\det (tE-A') \cdot \dot x=
\mbox{Ad}(tE-A')B' \cdot x= C \cdot x$ with polynomial matrix $C=(C_{i,j})$. From the structure
of $A'$ and $B'$, most of the entries in the last two columns of
$C$ are zeros, in particular: $C_{i,j}=0$, for $i=1, \ldots, n-2$ and $j=n-1,\
n$. This
means that the truncated curve $\overline \gamma (t)=(\oint_\delta \widetilde \omega_1, \ldots, \oint_\delta
\widetilde \omega_{n-2})$ satisfies the linear system whose matrix $\overline C$ is the
$(n-2) \times (n-2)$ upper-left corner of $C$. 

Starting from $q_0=(c_1, \ldots, c_{n-2}) \in \textbf C^{n-2}$, we derive the vectors $q_1, \ldots, q_{n-3} \in \textbf C[t]^{n-2}$ by: $q_{k+1}=Dq_k
+q_k \cdot \overline C$, with the same $D$ as before: $D= \det (tE-A')=(t-t_1)
\ldots (t-t_n)$. They satisfy: $D^k(q_0 \cdot \overline \gamma)=q_k \cdot \overline \gamma$.
 Let $\Delta$ be the wedge product of $q_0, \ldots,
q_{n-3}$. 

On the other hand, consider the matrix 
$$\widetilde {\mathcal P}=  
\left (
\begin{array}{ccc}
\oint_{\delta_1}  \widetilde \omega_1 & \cdots & \oint_{\delta_{n-2}}
\widetilde \omega_1\\
\vdots & \vdots & \vdots \\
\oint_{\delta_1} \widetilde \omega_{n-2} & \cdots & \oint_{\delta_{n-2}} \widetilde \omega_{n-2}
\end{array}
\right )$$

We obtain: 

\begin{equation}
\Delta \cdot \det \widetilde{\mathcal P}=  P^{\nu} \cdot \overline W
\end{equation}

with $P(t)=(t-t_1)\ldots (t-t_n)$, $\nu=(n-3)(n-2)/2$, and
$\overline W=W(\oint_{\delta_1}  \omega, \ldots, \oint_{\delta_{n-2}} \omega)$.

As $\overline W$ is nonzero (by the choice of the basis of the homology), both determinants $\Delta$ and
$\det  \widetilde{\mathcal P}$ are non identically vanishing. Notice that $\deg \Delta
\leq \frac{(n-1)(n-2)(n-3)}{2}$. A closer look at $\det \widetilde
{\mathcal P}$ shows that it is polynomial, of degree:

\begin{lem}

$\deg (\det  \widetilde {\mathcal P}) \leq n$.

\end{lem}

\begin{proof}
Note that the matrix $\widetilde {\mathcal P}$ is the product $R \cdot \overline{\mathcal P}$  of a $(n-2) \times
n$ constant matrix $R=(r_{i,j})$, $1 \leq i \leq n-2, 1\leq j \leq n$, of rank
  $n-2$  by $\overline{\mathcal P}$, obtained by removing the last two columns
  in the standard period matrix $\mathcal P$. There is no restriction in
  supposing the first $n-2$
  columns of $R$ independent (this amounts to permuting the cycles in $\overline{\mathcal
  P}$), and consider $\overline R$ the corresponding
  square matrix of rank $n-2$. Form the product $\overline R^{-1} \cdot  \widetilde {\mathcal
  P}$. Its determinant is the same as $\det \widetilde {\mathcal P}$, up to a
  nonzero constant. On the other hand, this matrix has the expression:  
$$\overline R^{-1} \cdot  \widetilde {\mathcal P}=  
\left (
\begin{array}{ccc}
\oint_{\delta_1} ( \omega_1 +\Omega_1) & \cdots & \oint_{\delta_{n-2}}
(\omega_1+ \Omega_1)\\
\vdots & \vdots & \vdots \\
\oint_{\delta_1} (\omega_{n-2} +\Omega_{n-2}) & \cdots & \oint_{\delta_{n-2}}
(\omega_{n-2}+ \Omega_{n-2})
\end{array}
\right )$$
where $\Omega_1, \ldots, \Omega_{n-2}$ belong to the span $\textbf C (\omega_{n-1},
\omega_n)$. Expanding $\det (\overline R^{-1} \cdot  \widetilde {\mathcal P})$,
it turns out that the term that brings the highest degree (with respect to
$t$) is the determinant:

$$  
\left \vert
\begin{array}{ccc}
\oint_{\delta_1}\Omega_1 & \cdots & \oint_{\delta_{n-2}} \Omega_1\\
\oint_{\delta_1} \omega_2 & \cdots & \oint_{\delta_{n-2}} \omega_2\\
\vdots & \vdots & \vdots \\
\oint_{\delta_1} \omega_{n-2} & \cdots & \oint_{\delta_{n-2}}
\omega_{n-2}
\end{array}
\right \vert$$ 

Setting
$x=t^{1/(n+1)} x'$, $y=t^{1/2} y'$, it follows that the leading term of this
determinant has degree $\frac{D}{n+1}$, where $D$ is the weighted degree 
$\deg \omega_n+ \deg \omega_2+
\ldots + \deg \omega_{n-2}= n+ \frac{n+1}{2}+\sum_{i=2}^{n-2}(i+\frac{n+1}{2})$,
hence $\frac{D}{n+1}\leq n$. The leading coefficient appears 
as the determinant of the integrals of the forms $\omega_n(x',y'),
\omega_2(x',y'), \ldots, \omega_{n-2}(x',y')$ over cycles in the level sets of the
\textit{principal} quasi-homogeneous part of $H$, $y^2-x^{n+1}$ (cf \cite{Gav2}). The latter
determinant is guaranteed to be nonzero since the differentials of the forms
involved 
are independent in the quotient ring of $\textbf C[x,y]$ by the Jacobian ideal
$(H_x, H_y)$.  
\end{proof}

\bigskip

Again, we observe that the identity (8) returns a
 \textit{quadratic} lower estimate: $ord_\infty \overline W \geq
\frac{n^2-7n+6}{2}$. On every finite pole: $ord_{t_i}\overline W
\geq 4-n$, and at a zero $t_0$ of $\overline W$: $ord_{t_0}\overline W\leq
\frac{n^2-n-6}{2}$. Finally: $ord_{t_0}\oint_\delta \omega \leq n-3+
 \frac{n^2-n-6}{2} \leq n-1 + \frac{n(n-1)}{2}$. The proof of
Theorem 2 is completed.

\bigskip

\subsection{Intersection with an algebraic hypersurface}

We consider the asymptotic behaviour of the integral with respect
to the degree $d$ of the form $\omega$.

\begin{thm}
Under the same assumptions on the Hamiltonian, consider the Abelian integral
of a 1-form $\omega$ of degree $d$. If $\oint_{\delta(t)} \omega \not \equiv 0$, then:
$$ord_{t_0} \left( \oint_{\delta(t)} \omega  \right)\leq A(n)+d \cdot B(n)$$
\end{thm}

\bigskip

An idea could be first to decompose $\omega$ in the Petrov module of $H$:
$\omega=p_1(t)\omega_1+ \ldots + p_n(t) \omega_n$, which implies:
$\oint_\delta \omega=p_1(t) \oint_\delta \omega_1+ \ldots + p_n(t)
 \oint_\delta \omega_n$, with $\deg p_i \leq \frac{d}{n+1}$,
 $i=1,\ldots,n$. Then apply the above reasoning, noticing that the vector
$$(\oint_{\delta} \omega_1, \ldots,
\oint_{\delta}
\omega_n, t\oint_{\delta} \omega_1, \ldots, t\oint_{\delta} \omega_n, \ldots,
t^{[d/(n+1)]} \oint_{\delta} \omega_1, \ldots,  t^{[d/(n+1)]} \oint_{\delta}
\omega_n)$$
 is a solution of a hypergeometric Picard-Fuchs system
of size $n \cdot ([d/(n+1)]+1)$.
 This would give a bound that is quadratic with
respect to the degree $d$ of the form. Yet, one should expect  linear growth,
since, for a fixed Hamiltonian, even the growth of the \textit{number} of zeros of the
integrals was proven by A. Khovanskii to be linear in the degree of the form. 

\bigskip

\begin{proof}

We define the curve $t \mapsto \Gamma(t)=(t, \oint_\delta \omega_1, \ldots,
\oint_\delta \omega_n) \subseteq \textbf C^{n+1}$, together with the algebraic
hypersurface $\{(t,x_1, \ldots, x_n) \in \textbf C^{n+1}$: $h(t, x_1, \ldots, x_n)=0\}$,
setting $h(t, x_1, \ldots, x_n)= p_1(t) x_1+ \ldots
+p_n(t) x_n$, in view of the above Petrov decomposition.
 Thus, $\oint_\delta \omega$ is the composition $(h \circ
\Gamma)(t)$. This reads also as the matrix product of the row vector
$Q_0=(p_1(t), \ldots, p_n(t)) \in \textbf C[t]^n$,
by the column vector $\gamma(t)=(\oint_{\delta} \omega_1, \ldots, 
\oint_{\delta}\omega_n)$.
The construction of vectors
$Q_k \in \textbf C [t]^n$ can be performed likewise. Their degrees, as well as the
degree of their exterior product, have affine growth with respect to $d$.
The lower estimate on the order of the Wronskian $W(\oint_{\delta_1} \omega, \ldots, 
\oint_{\delta_n}\omega)$ at its finite poles is not affected by $d$.   

\end{proof}

\end{document}